\documentclass[12pt]{article}
\usepackage{amsthm,amssymb,amsmath,txfonts,textcomp,mathrsfs,color}
\usepackage{mathrsfs}
\usepackage{amsmath}
\usepackage{amsfonts}
\usepackage{amssymb}
\usepackage{amsmath,amssymb,amsthm}
\usepackage{amsmath,amssymb,amsthm,amscd}
\usepackage[frame,cmtip,arrow,matrix,line,graph,curve]{xy}
\usepackage{graphpap, color}
\usepackage[mathscr]{eucal}
\usepackage{color}
\usepackage{verbatim}
\usepackage{hyperref}
\usepackage{cite}
\usepackage{slashbox,multirow}
\usepackage{float}
\newtheoremstyle{mythm}{3pt}{3pt}{}{}{\bfseries}{}{5mm}{}

\usepackage{float}
\catcode`!=11
\let\!iint\iint \def\iint{\displaystyle\!iint}
\let\!int\int   \def\int{\displaystyle\!int}
\let\!lim\lim   \def\lim{\displaystyle\!lim}
\let\!sum\sum   \def\sum{\displaystyle\!sum}
\let\!sup\sup   \def\sup{\displaystyle\!sup}
\let\!inf\inf   \def\inf{\displaystyle\!inf}
\let\!cap\cap   \def\cap{\displaystyle\!cap}
\let\!max\max   \def\max{\displaystyle\!max}
\catcode`!=12

\newtheorem{thm}{Theorem}[section]
\newtheorem{lem}[thm]{Lemma}

\newtheorem*{ty}{ Theorem (Tang and Yan)}

\theoremstyle{remark}\newtheorem{rem}{Remark}[section]
\newtheorem*{ack}{Acknowledgment}

\allowdisplaybreaks

\numberwithin{equation}{section}
\begin{document}

\title{The First Eigenvalue for the Bi-Beltrami-Laplacian \\ on Minimal Isoparametric Hypersurfaces of $\mathbb{S}^{n+1}(1)$}
\author{Lingzhong Zeng
\\ \small   College of Mathematics and Informational Science,
Jiangxi Normal University,
\\
\small Nanchang 330022, China, E-Mail: lingzhongzeng@yeah.net
}
\date{}

\maketitle

\begin{abstract}\noindent In this paper, we investigate the first eigenvalues of two closed eigenvalue problems of the bi-Beltrami-Laplacian on minimal embedded isoparametric hypersurface in the unit sphere
$\mathbb{S}^{n+1}(1)$. Although many mathematicians want to derive the corresponding results for the first eigenvalues of bi-Beltrami-Laplacian, they encountered great difficulties in proving the limit theorem of the version of bi-Beltrami-Laplacian along with the strategy due to I. Chavel and E. A. Feldman(Journal of Functional Analysis, 30 (1978), 198-222) and S. Ozawa (Duke Mathematics  Journal, 48 (1981),767-778). Therefore, as the author knows, there are no any results of Tang-Yan type( Journal of Differential Geometry, 94 (2013) 521-540).  However, by  the variational argument, we overcome the difficulties and determine the first eigenvalues of the bi-Beltrami-Laplacian in the sense of  isoparametric hypersurfaces.  We note that our proof is quite simple.


\vskip3mm
\noindent 2000 MSC 35P15, 53C40.
\end{abstract}

\section{Introduction}
Let $(M^{n},g)$ be an $n$-dimensional complete  Riemannian
manifolds. For any $C^{k},k\geq2$ function $f$ on $M^{n}$,  define the
Beltrami-Laplacian of $f$ by
$$\Delta=-\frac{1}{\sqrt{det(g)}}\sum_{i,j=1}\partial_{i}g^{ji}\sqrt{det(g)}\partial_{j},$$and then, it is a positive elliptic operator. We consider the following closed eigenvalue problem:
\begin{equation}\label{L} \Delta u=\lambda u.
\end{equation}
When $M^{n}$ is a minimal embedded hypersurface in the unit sphere
$\mathbb{S}^{n+1}(1)$, it follows from Takahashi Theorem that
$\lambda_{1}(M)$ is not greater than $n$. In this connection,
S.-T. Yau proposed in 1982 the following conjecture{\rm\cite{SY,Yau}}:
The first eigenvalue of every closed minimal
hypersurface $M^{n}$ in the unit sphere $S^{n+1}(1)$ is just the dimension of the hypersurface  $M^{n}$.
Attacking the Yau's conjecture, many mathematicians contributed to this problem, we refer to  \cite{CW,SY,MOU,Kot,Sol1, Sol2,TY2} and reference therein. In particular, based on the generic far reaching results on the classification of isoparametric
hypersurfaces in $\mathbb{S}^{n+1}(1)$ (we refer to \cite{CCJ,Imm,Chi,DN,Mu1,Mu2,OT1,OT2,Miy} and the references therein), Tang and Yan \cite{TY2} made a landmark breakthrough to Yau's conjecture.
They considered a little more restricted problem of Yau's conjecture for
closed minimal isoparametric hypersurfaces $M^{n}$ in
$\mathbb{S}^{n+1}(1)$.  In detail, they proved the following:

\begin{ty}Let $M^{n}$ be a closed minimal isoparametric hypersurface in the
unit sphere $\mathbb{S}^{n+1}(1)$. Then, $\lambda_{1}(M^{n}) =
n$.\end{ty}Up to now,
Yau's conjecture is still open.

\vskip 2mm
Next, we consider the following two closed eigenvalue problems of bi-Beltrami-Laplaican:
 \begin{equation}\label{C}\Delta^{2}u=\Lambda u;
\end{equation}and
\begin{equation}\label{B} \Delta^{2}u=\Gamma\Delta u.
\end{equation}
Let $\Lambda_{i}$ and $\Gamma_{i}$ denote the $i$-th eigenvalue of eigenvalue problem \eqref{C} and eigenvalue problem \eqref{B}, respectively. Then, the eigenvalues of this eigenvalue problems \eqref{C} and \eqref{B} are real and discrete:
$$0=\Lambda_{0} \leq \Lambda_{1}\leq\Lambda_{2}\leq\Lambda_{3}\leq\cdots\rightarrow+\infty,$$ and$$0=\Gamma_{0} \leq \Gamma_{1}\leq\Gamma_{2}\leq\Gamma_{3}\leq\cdots\rightarrow+\infty,$$
where each $\Lambda_{i}$ and each $\Gamma_{i}$ have finite multiplicity.
In \cite{TY2}, Tang and Yan noted that\emph{ the calculation of the spectrum of the Laplace-Beltrami operator, even of the first eigenvalue, is rather complicated
and difficult.} In fact, the author  thinks that, comparing with the case of the Laplace-Beltrami operator, the calculation of the first eigenvalue of the bi-Laplace-Beltrami operator maybe more  complicated and difficult. In order to derive the results for the first eigenvalues of the bi-Beltrami-Laplacian on isoparametric hypersurface which is embedded into the $(n+1)$-dimensional unit sphere $\mathbb{S}^{n+1}$, the natural way is to use a series of techniques due to Muto-Ohnita-Urakawa \cite{MOU}, Kotani
\cite{Kot}, and Solomon \cite{Sol1, Sol2}, Muto \cite{Mut} and Tang-Yan \cite{TY2}. Therefore, firstly, it seems unavoidable to prove a crucial limit theorem in the sense of bi-Beltrami-Laplacian along with the strategy due to I. Chavel and E. A. Feldman \cite{CF} and S. Ozawa\cite{Oza}. In fact, many mathematicians considered the problem, but, technically, they encountered great difficulty in proving the limit theorem. Therefore, as we know,
there are no any results of Tang-Yan type in the sense of bi-Laplace-Beltrami operator. However, by an observation of variational principle, we find that it is not necessary to prove the limit theorem but provide an alternative strategy to overcome those difficulties.  More precisely, by variational principle and based on Tang-Yan's work\cite{TY2}, we prove that following:
\begin{thm}\label{thm-z}Let $M^{n}$ be a closed minimal isoparametric hypersurface in the
unit sphere $\mathbb{S}^{n+1}(1)$. Then $$\Lambda_{1}(M^{n}) =
n^{2},\ \ and \ \ \Gamma_{1}(M^{n}) =
n.$$\end{thm}

\section{The Proof of Main Result}
In this section, we give a very short proof of theorem \ref{thm-z}. Firstly, we need the following theorem \cite{T}:
\begin{thm}\emph{\textbf{(Takahashi)}} Let $M^{n}$ be an $n$-dimensional Riemannian manifold. For an isometric immersion
$\psi:M^{n}\rightarrow\mathbb{S}^{n+1}$ it is minimal immersion into
$\mathbb{S}^{n}$ if and only if
\begin{equation}\label{t1}\Delta\psi=n\psi.\end{equation}
\end{thm}By the variational principle and Schwarz's inequality, it is easy to prove the following lemma:
\begin{lem}\label{lem}Let $\lambda_{i}$, $\Lambda_{i}$ and $\Gamma_{i}$ denote the $i$-th eigenvalues of eigenvalue problem \eqref{L}, eigenvalue problem \eqref{C}  and eigenvalue problem \eqref{B}, respectively. Then,\begin{equation}\label{v1} \Lambda_{k}\geq\lambda_{k}^{2},\ \
  and \ \ \Gamma_{k}\geq\lambda_{k}.\end{equation}
\end{lem}
\textbf{\emph{ Proof of theorem {\rm\ref{thm-z}.}} }
In \eqref{t1}, for any $i,(i=1,2,\cdots,n+2)$, we take the standard coordinate function $\psi=x_{i}$, and then, one has  $\Delta x_{i}=n x_{i}$, and thus, $\Delta^{2} x_{i}=n\Delta x_{i}$. Using Takahashi's theorem again, it follows that
$\Delta^{2} x_{i}=n^{2} x_{i}$.
Therefore, by the variational principle, we have
$\Lambda_{1}\leq n^{2}$ and $\Gamma_{1}\leq n$. Furthermore, the inequalities \eqref{v1} in lemma \ref{lem} implies  $\Lambda_{1}\geq n^{2}$ and $\Gamma_{1}\geq n$. Hence, we have $ \Lambda_{1}= n^{2}$ and $\Gamma_{1}= n$.  The proof of theorem \ref{thm-z} is complete.
$$\eqno\Box$$

\begin{rem}By the observation for the proof, it is not difficult to see that we can remove the condition of isoparametric hypersurfaces in theorem \ref{thm-z}, this is, both $\Lambda_{1}(M^{n})=n^{2}$ and $ \Gamma_{1}(M^{n})=n$ hold without the assumption of the isoparametric  condition, if one can prove that Yau's conjecture is ture. \end{rem}

\begin{rem}According to the variational principle and Choi and Wang's result \cite{CW}, we can easily prove that, $$\Lambda_{1}(M^{n}) \geq
\frac{n^{2}}{4},\ \ and \ \ \Gamma_{1}(M^{n})\geq
\frac{n}{2},$$ without the assumption of isoparametric hypersurfaces. \end{rem}

\begin{rem}For eigenvalues of the closed eigenvalue problem \eqref{C} and closed eigenvalue problem \eqref{B} of the Bi-Laplacian on
focal submanifolds of isoparametric hypersurfaces in the unit sphere, we can obtain
similar results by making use of the variational principle and the results on the first eigenvalue of Laplacian $\lambda_{1}$ due to Tang and
Yan \cite{TY2} and Tang, Xie and Yan \cite{TXY}.\end{rem}

\begin{ack} The author is supported by the National Nature Science Foundation of China (Grant No. 11401268).\end{ack}

\end{document}